\documentclass[11pt, final]{article}
\usepackage{a4}
\usepackage{amsmath}
\usepackage{amstext}
\usepackage{amssymb}
\usepackage{showkeys}
\usepackage{epsfig}
\usepackage{cite}
\newtheorem{theorem}{Theorem}
\newtheorem{corollary}[theorem]{Corollary}
\newtheorem{df}{Definition}
\newenvironment{definition}{\df\rm}{\enddf}
\newtheorem{ex}{Example}
\newenvironment{example}{\ex\rm}{\endex}
\newtheorem{lemma}[theorem]{Lemma}
\newtheorem{proposition}[theorem]{Proposition}
\newtheorem{rem}{Remark}
\newenvironment{remark}{\rem\rm}{\endrem}

\newcounter{unnumber}

\newtheorem{prf}[unnumber]{Proof.}
\newenvironment{proof}{\prf\rm}{\hfill{$\blacksquare$}\endprf}
\newcommand{\R}{\mathbb{R}}%
\newcommand{\N}{\mathbb{N}}%
\newcommand{\ra}{\rightarrow}%

\DeclareMathOperator*\inte{int}%
\DeclareMathOperator*\cl{cl}%
\DeclareMathOperator*\gap{gap}%
\DeclareMathOperator*\epi{epi}%
\DeclareMathOperator*\dom{dom}%
\DeclareMathOperator*\Limsup{Limsup}%

\DeclareMathOperator*\lm{\longrightarrow}%

\title{On the Dini-Hadamard subdifferential of the difference of two functions}

\author{Radu Ioan Bo\c{t} \thanks
{Faculty of Mathematics, Chemnitz University of Technology,
D-09107 Chemnitz, Germany, e-mail:
radu.bot@mathematik.tu-chemnitz.de. Research partially supported by DFG (German Research Foundation), project WA 922/1-3.} \and Delia-Maria Nechita
\thanks {Faculty of Mathematics and Computer Science, Babe\c{s}-Bolyai University, Cluj-Napoca,
Romania, e-mail: delia-maria.nechita@s2009.tu-chemnitz.de. Research done during the stay of the author in the academic year 2009/2010 at Chemnitz University of Technology
as a guest of the Chair of Applied Mathematics (Approximation Theory). The author wishes to thank for the financial support provided from programs co-financed by The Sectoral Operational Programme Human Resources Development,
Contract POSDRU 6/1.5/S/3 –- ``Doctoral studies: through science towards society''.}}

\begin{document}
\maketitle

\noindent \textbf{Abstract.} In this paper we first provide a general formula of inclusion for the Dini-Hadamard $\varepsilon$-subdifferential of the difference of two functions and show that it becomes equality in case the functions are directionally approximately starshaped at a given point and a weak topological assumption is fulfilled. To this end we give a useful characterization of the Dini-Hadamard $\varepsilon$-subdifferential by means of sponges. The achieved results are employed in the formulation of optimality conditions via the Dini-Hadamard subdifferential for cone-constrained optimization problems having the difference of two functions as objective.\\



\section{Introduction}

In this paper we provide a formula for the \textit{Dini-Hadamard $\varepsilon$-subdifferential} of the difference of two functions by making use of the star-difference of the Dini-Hadamard $\varepsilon$-subdifferentials of the functions involved.  In this investigation an important role will be played by a variational description of the \textit{Dini-Hadamard $\varepsilon$-subgradients} of an arbitrary function via \textit{sponges}, a notion introduced by Treiman in \cite{Treiman}, which represents the counterpart of a well-known variational description for \textit{Fr\' echet $\varepsilon$-subgradients}.

While in the announced subdifferential formula for the difference of two functions one inclusion follows automatically, in order to guarantee the other one we need some supplementary assumptions on the functions involved. More precisely, we show that in case the two functions are directionally approximately starshaped at a given point and a weak topological assumption is fulfilled, then the opposite inclusion is fulfilled, too. The class of \textit{directionally approximately starshaped functions} has been introduced in \cite{Penot-as-das} and contains the class of \textit{approximately starshaped} and, consequently, the one of \textit{approximately convex functions} (see \cite{Aussel, Ngai-ac, Penot-ac} for more on these classes of functions). We also give a characterization of the directionally approximately starshaped functions by means of sponges and furnish an example of a function which is directionally approximately starshaped at a given point, but not approximately starshaped at that point. We also show that the Fr\' echet subdifferential and the Dini-Hadamard subdifferential of this function at the point in discussion do not coincide. The weak topological assumption we use in the main result is an extension of the \textit{gap-continuity} introduced by Penot in \cite{Penot-gap-continuity}. The results in this article extend some assertions given in \cite{Amahroq} for approximately starshaped functions via the Fr\' echet subdifferential and are motivated by the fact that there exist directionally approximately starshaped functions for which the statements in \cite{Amahroq} do not apply.

Finally, we employ these results in the formulation of optimality conditions for cone-constrained optimization problems having the difference of two functions as objective and give necessary and sufficient conditions for the existence of so-called \textit{spongiously local $\varepsilon$-blunt minimizers} for all $\varepsilon > 0$, a notion which represents an extension of the \textit{local $\varepsilon$-blunt minimizer} introduced and investigated in \cite{Amahroq}.

The paper is organized as follows. In Section 2 we introduce some preliminary notions and results, give the variational description of the Dini-Hadamard $\varepsilon$-subdifferential by means of sponges and provide a partial result for the Dini-Hadamard $\varepsilon$-subdifferential of the difference of two functions. In the following section we first introduce the class of directionally approximately starshaped functions and study its relations with the ones of
approximately starshaped and approximately convex functions. Afterwards we prove the main result of the paper, the characterization of the directionally approximately starshaped functions via sponges playing here a determinant role. Finally, we turn our attention to the formulation of optimality conditions via the Dini-Hadamard subdifferential for a general cone-constrained optimization problem having the difference of two functions as objective.

\section{Preliminaries}

Consider a Banach space $X$ and its topological dual space $X^*$. We denote the \textit{open ball} with center $\overline{x} \in X$ and radius $\delta >0$ in $X$ by $B(\overline{x},\delta)$, while $\overline{B}_X$ and $S_X$ stand for the \textit{closed unit ball} and the \textit{unit sphere} of $X$, respectively. Having a set $C \subseteq X$, $\delta_C : X \rightarrow  \R \cup \{+\infty\}$, defined by $\delta_C(x) = 0$ for $x \in C$ and $\delta_C(x) = +\infty$, otherwise, denotes its \textit{indicator function}.

Let $f:X\ra\R\cup\{+\infty\}$ be a given function. As usual, we denote by $\dom f = \{ x \in X: \ f(x)<+ \infty \}$ the \emph{effective domain} of $f$ and by $\epi f
=\{(x,\alpha) \in X \times \R: \ f(x)\leq \alpha \}$  the \emph{epigraph} of $f$. We make the conventions $(+\infty)-(+\infty)=+\infty$ and $0(+\infty)=+\infty$. For $\varepsilon \geq 0$ the \emph{Fr\'{e}chet $\varepsilon$-subdifferential} (or the \emph{analytic $\varepsilon$-subdifferential}) of $f$ at $\overline{x} \in \dom{f}$ is defined by
$$\partial^F_{\varepsilon} f(\overline{x}):= \left\{ x^* \in X^* : \liminf_{\|h\| \ra 0} \frac{f(\overline{x}+h)-f(\overline{x})-\langle x^*,h \rangle}{\|h\|} \geq - \varepsilon \right \},$$
which means that one has
\begin{equation}\label{charFr}
\begin{array}{rl}
\overline{x}^* \in \partial^F_{\varepsilon} f(\overline{x}) \Leftrightarrow & \ \mbox{for all} \ \alpha >0 \ \mbox{there exists} \ \delta >0 \ \mbox{such that} \ \mbox{for all} \ x \in B(\overline{x}, \delta)\\
& f(x)-f(\overline{x}) \geq \langle \overline{x}^*, x-\overline{x} \rangle - (\alpha + \varepsilon) \|x-\overline{x}\|.
\end{array}
\end{equation}

The following constructions
$$d^-f(\overline{x};h):=\liminf_{\substack{ u \ra h \\ t \downarrow 0}} \frac{f(\overline{x}+tu)-f(\overline{x})} {t}=\sup_{\substack{\delta >0}} \inf_{\substack{u \in B(h, \delta) \\ t \in (0, \delta)}} \frac{f(\overline{x}+tu)-f(\overline{x})}{t}$$
and (see \cite{Ioffe1, Ioffe2})
$$\partial^-_{\varepsilon} f(\overline{x}):=\{ x^* \in X^* : \langle x^*,h \rangle \leq d^-f(\overline{x};h) + \varepsilon \|h\| \ \mbox{for all} \ h \in X\}, \ \mbox{where} \ \varepsilon \geq 0,$$
are called the \emph{Dini-Hadamard directional derivative} of $f$ at $\overline{x}$ in the direction $h \in X$ and the \emph{Dini-Hadamard $\varepsilon$-subdifferen\-tial} of $f$ at $\overline{x}$, respectively. When $\overline{x} \not \in \dom f$ we set $\partial^F_{\varepsilon} f(\overline{x})=\partial^-_{\varepsilon} f(\overline{x}):=\emptyset$ for all $\varepsilon \geq 0$. Note that for $\overline{x} \in \dom{f}$ the function $d^-f(\overline{x};\cdot)$ is in general not convex, while $\partial^-_{\varepsilon} f(\overline{x})$ is always a convex set. Further, we notice that $d^-f(\overline{x};0)$ is either $0$ or $-\infty$ (see \cite{Ioffe-df}).

The function $f:X\ra\R\cup\{+\infty\}$ is said to be \emph{calm at $\overline{x} \in \dom f$} if there exists $c \geq 0$ and $\delta > 0$ such that
$f(x) - f(\overline{x}) \geq -c\|x-\overline{x}\|$ for all $x \in B(\overline{x}, \delta)$. For $\overline{x} \in \dom f$ one has (see, for instance, \cite{giner}) that $f$ is calm at $\overline{x}$ if and only if $d^-f(\overline{x};0) = 0$.

For $\varepsilon=0$, $\partial^-f(\overline{x}):=\partial^-_{0} f(\overline{x})$ and $\partial^F f(\overline{x}):= \partial^F_{0} f(\overline{x})$ are nothing else than the \emph{Dini-Hadamard subdifferential} and the \emph{Fr\'{e}chet subdifferential} of $f$ at $\overline{x}$, respectively.

For all $\varepsilon \geq 0$ we have the following inclusion $$\partial^F_{\varepsilon} f(\overline{x}) \subseteq \partial^-_{\varepsilon} f(\overline{x}),$$
which may be in general strict.

To see this in case $\varepsilon = 0$ consider, for instance, the function $f: C[0,1] \ra \R, \ f(x)= -\|x\|_\infty$. One has that
$\partial^F f(x) = \emptyset$ for all $x \in C[0,1]$, while $\partial^-f(\overline{x}) \neq \emptyset$ when $\overline{x} \in S_{C[0,1]}$ is chosen such that $|\overline{x}| : [0,1] \rightarrow \R, |\overline x|(t) = |\overline{x}(t)|$, attains its maximum at exactly one point of the interval $[0,1]$ (see \cite[Exercise 8.28]{Fabian}). For a similar example in $\ell_1$ we refer to \cite[Exercise 8.26]{Fabian}. However, when $X$ is finite dimensional, we always have that $\partial^F f(x)=\partial^- f(x)$ for all $x \in X$.

The next subdifferential notion we need to recall is the one of \emph{$G$-subdifferential} and we describe in the following the procedure of constructing it (see \cite{Ioffe3}). To this aim we consider first the \emph{$A$-subdifferential} of $f$ at $\overline{x} \in \dom f$, which is defined via topological limits as follows
$$\partial^A f(\overline{x}):=\bigcap_{\substack{L \in \mathcal{F}(X)\\ \varepsilon > 0}} \overline{\Limsup_{\substack{x \lm\limits^{f} \overline{x}}}} \ \partial^-_{\varepsilon} (f + \delta_{x+L})(x),$$
where $\mathcal{F}(X)$ denotes the collection of all finite dimensional subspaces of $X$ and $\overline{\Limsup}$ stands for the \textit{topological counterpart} of the \textit{sequential Painlev\'{e}-Kuratowski upper/outer limit} of a set-valued mapping with sequences replaced by nets and where $x \lm\limits^{f} \overline{x}$ means $x \lm \overline x$ and $f(x) \lm f(\overline x)$. More precisely, for a multivalued mapping $F: X \rightrightarrows X^*$, we say that $x^* \in \overline{\Limsup}_{\substack{x \ra \overline{x}}} \ F(x)$ if for each weak$^*$-neighborhood $\mathcal{U}$ of the origin of $X^*$ and for each $\delta > 0$ there exists $x \in B(\overline x, \delta)$ such that $(x^* + \mathcal{U})\cap F(x) \neq \emptyset$.

The \emph{$G$-normal cone} to a set $C \subseteq X$ at $\overline{x} \in C$ is defined as
$$N^G(C,\overline{x}):=\cl\nolimits^*\left(\bigcup_{\substack{\lambda >0}} \lambda \partial^A d(\overline{x},C)\right),$$
where $d(\overline{x},C):=\inf_{c \in C}\|\overline x - c\|$ denotes the \textit{distance} from $\overline{x}$ to $C$ and $\cl^*$ stands for the \textit{weak$^*$-closure} of a set in $X^*$, while the \emph{$G$-subdifferential} of $f$ at $\overline{x} \in \dom f$ can be defined now as follows
$$\partial^G f(\overline{x}):=\left\{ x^* \in X^*: \ (x^*, -1) \in N^G(\epi f, (\overline{x},f(\overline{x}))) \right\}.$$
When $\overline{x} \not \in \dom f$ we set $\partial^A f(\overline{x})=\partial^G f(\overline{x}):=\emptyset$. Thus, by taking into account \cite[Proposition 4.2]{Ioffe3}, we have for all $x \in X$ the following relations of inclusion
\begin{equation}\label{FDG} \partial^F f(x) \subseteq \partial^- f(x) \subseteq \partial^G f(x).\end{equation}

One can notice that when $f$ is a convex function it holds $\partial^F f(x)=\partial^- f(x)=\partial^G f(x)=\partial f(x)$ for all $x \in X$, where $\partial f (\overline x):=\{x^*\in X^*: f(y)-f(\overline x)\geq \langle x^*, y- \overline x \rangle \ \forall y\in X\}$,  for $\overline x \in \dom f$, and $\partial f (\overline x):=\emptyset$, otherwise, denotes the subdifferential of $f$ at $\overline x$ in the sense of convex analysis.

For $f: X \ra \R \cup \{+\infty\},$ $\overline{x} \in \dom f$ and $\varepsilon \geq 0$, we define $f_\varepsilon: X \ra \R \cup \{+\infty\}$ as being
\begin{equation}\label{f epsilon} f_\varepsilon (x):=f(x)+\varepsilon \|x-\overline{x}\|. \end{equation}

\begin{lemma}\label{dini fepsilon} Let $f:X \ra \R \cup \{+\infty\}$ be a given function and $\overline{x} \in \dom f$. Then for all $\varepsilon \geq 0$ it holds
\begin{equation}\label{fepsilon dini} \partial^-_{\varepsilon} f(\overline{x})=\partial^- f_\epsilon(\overline{x}). \end{equation}
\end{lemma}

\begin{proof} Because of the equivalence
\begin{eqnarray*}
    x^* \in \partial ^- f_\varepsilon(\overline{x})
    \Leftrightarrow \mbox{for all} \ h \in X: \ \langle x^*,h \rangle \leq d^- f_{\varepsilon} (\overline{x},h),
\end{eqnarray*}
it is sufficient to prove only the equality $d^- f_{\varepsilon} (\overline{x},h)=d^- f (\overline{x},h) + \epsilon \|h\|$ for arbitrary $h \in X$. Let $h \in X$. Then we have
\begin{eqnarray*}
    d^- f_{\epsilon} (\overline{x},h) &=& \liminf_{\substack{u \ra h \\ t \downarrow 0}} \frac{f_\varepsilon(\overline{x} + tu) - f_\epsilon (\overline{x})}{t} \\
    &=& \liminf_{\substack{u \ra h \\ t \downarrow 0}} \left(\frac{f(\overline{x} + tu) - f(\overline{x})}{t} + \epsilon \|u\|\right) \\
    &=& \liminf_{\substack{u \ra h \\ t \downarrow 0}} \frac{f(\overline{x} + tu) - f(\overline{x})}{t} + \lim_{\substack{u \rightarrow h \\ t \downarrow 0}} \epsilon \|u\| \\
    &=& d^- f (\overline{x},h) + \epsilon \|h\|,
\end{eqnarray*}
which completes the proof.
\end{proof}

\begin{remark}\label{remarca-1} Let us notice that one can replace in \eqref{fepsilon dini} the Dini-Hadamard subdifferential by the Fr\'{e}chet one (see, for instance, \cite{Amahroq}). Thus, when $X$ is finite dimensional, $\partial^F_{\varepsilon} f(x)=\partial^-_{\varepsilon} f(x)$ for all $x \in X$ and all $\varepsilon \geq 0$. On the other hand, in case $f$ is convex, by a classical subdifferential sum formula provided by the convex analysis, one has
$\partial^- f_{\varepsilon}(x)=\partial f(x)+ \varepsilon \overline{B}_{X^*}$ for all $x \in X$ and all $\varepsilon \geq 0$. Finally, we notice that for the Dini-Hadamard $\varepsilon$-subdifferential the following monotonicity property holds:
\begin{eqnarray}\label{monotonicity}
\partial^-_{\varepsilon_1} f(x) \subseteq \partial^-_{\varepsilon_2} f(x),
\end{eqnarray}
when $\varepsilon_2 \geq \varepsilon_1 \geq 0$ and $x \in X$.
\end{remark}

The notion that we recall below was introduced by Treiman in \cite{Treiman} and, as we will prove in Theorem \ref{sponge-subdifferential}, it turns out to be very useful for characterizing the Dini-Hadamard subdifferential.

\begin{definition}\label{sponge-def}
\rm A set $S \subseteq X$ is said to be a \emph{sponge around $\overline{x} \in X$} if for all $h \in X \setminus \{0\}$ there exist $\lambda >0$ and $\delta >0$ such that $\overline{x}+[0,\lambda] \cdot B(h,\delta) \subseteq S$.
\end{definition}

\begin{example}\label{ex-sponge} (cf. \cite[Example 2.3]{Treiman})
\rm Let $f:X\ra \R$ be a locally Lipschitz and G\^ateaux differentiable function at $\overline{x} \in X$ with $x^* \in X^*$ its G\^ateaux derivative at this point. Then for all $\varepsilon >0$ the sets
$$S_1:=\{x \in X: \ f(x)-f(\overline{x}) \geq \langle x^*, x-\overline{x} \rangle - \varepsilon \|x-\overline{x}\|\}$$
and
$$S_2:=\{x \in X: \ f(x)-f(\overline{x}) \leq \langle x^*, x-\overline{x} \rangle + \varepsilon \|x-\overline{x}\|\}$$
are sponges around $\overline{x}$.
\end{example}

\begin{remark}\label{remarca-sponge}
\rm Every neighborhood of a point $\overline{x} \in X$ is also a sponge around $\overline{x}$, but the converse is not true (see for instance Example \ref{sponge-not-neighborhood} below). Nevertheless, when $S$ is a convex set or $X$ is a finite dimensional space (here one can make use of the fact that the unit sphere is compact), then every sponge around $\overline{x}$ is also a neighborhood of $\overline{x}$.
\end{remark}

The following notion has been introduced in \cite{Treiman}, too.

\begin{definition}\label{he-subgradient}
Let $f:X \ra \R \cup \{+\infty\}$ be a given function, $\overline x \in \dom f$ and $\varepsilon \geq 0$. We say that $x^* \in X^*$ is an \emph{$H_\varepsilon$-subgradient} of $f$ at $\overline{x}$ if there exists a sponge $S$ around $\overline{x}$ such that for all $x \in S$
\begin{eqnarray*}
f(x)-f(\overline{x}) \geq \langle x^*, x-\overline{x} \rangle - \varepsilon \|x-\overline{x}\|.
\end{eqnarray*}
\end{definition}

The following lemma was inspired by some statements one can find in Treiman's paper \cite{Treiman}.

\begin{lemma}\label{treiman}
Let $f:X \ra \R \cup \{+\infty\}$ be a given function, $\overline x \in \dom f$ and $\varepsilon \geq 0$. The following statements are
true:

(i) If $x^* \in \partial^-_{\varepsilon} f(\overline{x})$, then $x^*$ is an $H_\gamma$-subgradient of $f$ at $\overline{x}$ for all $\gamma > \varepsilon$.

(ii) If $f$ is calm at $\overline{x}$ and $x^*$ is an $H_\varepsilon$-subgradient of $f$ at $\overline{x}$, then $x^* \in \partial^-_{\varepsilon} f(\overline{x})$.
\end{lemma}

\begin{proof}(i) Let $x^* \in \partial^-_{\varepsilon} f(\overline{x})$ and $\gamma > \varepsilon$ be given. We consider the set
$$S:=\{ x \in X: \ f(x)-f(\overline{x}) \geq \langle x^*, x-\overline{x}  \rangle -\gamma \|x-\overline{x}\| \}$$
and show in the following that $S$ is a sponge around $\overline{x}$, which will complete the proof of the first statement. To this end, we fix an arbitrary element $h \in X \setminus \{0\}$. Since $x^* \in \partial^-_{\varepsilon} f(\overline{x})$, we obtain
\begin{eqnarray*}
    \liminf_{\substack{u \ra h \\ t \downarrow 0}} \frac{f(\overline{x} + tu) - f(\overline{x})}{t} &\geq& \langle x^*, h  \rangle - \varepsilon \|h\| \\
    &>& \langle x^*, h  \rangle - \left(\frac{\gamma + \varepsilon}{2}\right) \|h\| > \langle x^*, h  \rangle - \gamma \|h\|.
\end{eqnarray*}
Thus there exists $\delta_1 >0$ such that for all $u \in B(h,\delta_1)$
$$\langle x^*, h  \rangle - \left(\frac{\gamma + \varepsilon}{2}\right) \|h\| > \langle x^*, u  \rangle - \gamma \|u\|.$$
On the other hand, using the definition of the limit inferior, there exist $\delta_2>0$ such that for all $t \in (0,\delta_2)$ and all $u \in B(h,\delta_2)$ it holds
$$\frac{f(\overline{x} + tu) - f(\overline{x})}{t} > \langle x^*, h  \rangle - \left(\frac{\gamma + \varepsilon}{2}\right) \|h\|.$$
Hence, the last two relations above lead to the existence of $\delta:= \frac{1}{2}\min \{\delta_1, \delta_2 \}>0$ such that for all $t \in (0,\delta]$ and all $u \in B(h, \delta)$
$$f(\overline{x}+tu)-f(\overline{x})>\langle x^*, tu  \rangle - \gamma \|tu\|.$$
Now, it is not difficult to see that for all $x\in \overline{x}+[0,\delta]B(h,\delta)$ one gets
$$f(x)-f(\overline{x}) \geq \langle x^*, x-\overline{x}  \rangle - \gamma \|x-\overline{x}\|$$
and therefore $\overline{x}+[0,\delta]B(h,\delta) \in S$, which means in fact that $S$ is a sponge around $\overline{x}$.

(ii) Take $x^*$ to be an arbitrary $H_\varepsilon$-subgradient of $f$ at $\overline{x}$. Via Definition \ref{he-subgradient} one gets the existence of a sponge $S$ around $\overline{x}$ such that for all $x \in S$
$$f(x)-f(\overline{x}) \geq \langle x^*, x-\overline{x} \rangle - \varepsilon \|x-\overline{x}\|.$$
Let $h \in X \setminus \{0\}$ be fixed. Then there exist $\lambda>0$ and $\delta >0$ such that for all $t \in (0, \lambda]$ and $u \in B(h, \delta)$, one has $\overline{x}+tu \in S$ and, consequently,
$$\frac{f(\overline{x}+tu)-f(\overline{x})}{t} \geq \langle x^*,u \rangle - \varepsilon \|u\|.$$
Passing now to limit inferior in both sides, we obtain
$$\liminf_{\substack{u \ra h \\ t \downarrow 0}} \frac{f(\overline{x} + tu) - f(\overline{x})}{t} \geq \liminf_{\substack{u \ra h \\ t \downarrow 0}} \left[ \langle x^*,u \rangle -\varepsilon \|u\| \right]=\langle x^*,h \rangle -\varepsilon \|h\|.$$
This actually means that for all $h \in X \setminus \{0\}$
\begin{eqnarray}\label{h-esubdifferential}
d^-f(\overline{x};h) \geq \langle x^*,h \rangle - \varepsilon \|h\|.
\end{eqnarray}
As $f$ is calm at $\overline{x}$, it holds that $d^-f(\overline{x};0) = 0$ and, consequently, $x^* \in  \partial^-_{\varepsilon} f(\overline{x})$.
\end{proof}

\begin{remark}\label{sponge}
From the proof of Lemma \ref{treiman} (i) one can easily conclude that whenever $\overline x \in \dom f$, $\varepsilon \geq 0$, $x^* \in \partial_\varepsilon^-f(\overline{x})$ and $\gamma > \varepsilon$ the set \begin{eqnarray}\label{eq7}
S:=\{x \in X: \ f(x)-f(\overline{x}) \geq \langle x^*, x-\overline{x}\rangle -\gamma \|x-\overline{x}\| \}
\end{eqnarray}
is a sponge around $\overline{x}$.
\end{remark}

In the following we provide an example of a sponge around a point which is not a neighborhood of that point.

\begin{example}\label{sponge-not-neighborhood}
We consider again the space $C[0,1]$ endowed with the supremum norm. Let $\overline{x} \in S_{C[0,1]}$ be an element in this space with the property that $|\overline{x}|$ attains its maximum at exactly one point of the interval $[0,1]$. Let further $x^* \in X^*$ be an element in $\partial^-(-\|\cdot\|_\infty) (\overline{x})$, which, as we have seen, is a nonempty set. As the Fr\'echet subdifferential of $-\|\cdot\|_\infty$ at $\overline{x}$ is empty, there exists an $\alpha >0$ such that for all $\delta>0$ there is some $x \in B(\overline{x}, \delta)$ satisfying
\begin{eqnarray}\label{s}
\|\overline{x}\|_\infty-\|x\|_\infty + \alpha \|x-\overline{x}\|_\infty< \langle x^*, x-\overline{x}\rangle.
\end{eqnarray}
As seen above, the set
\begin{eqnarray}
S:=\{x \in C[0,1]: \ \|\overline{x}\|_\infty-\|x\|_\infty+\alpha\|x-\overline{x}\|_\infty \geq \langle x^*, x-\overline{x}\rangle \}
\end{eqnarray}
is a sponge around $\overline{x}$ (take $\varepsilon :=0$ and $\gamma:=\alpha$ in Remark \ref{sponge}). It remains to show that  $S$ is not a neighborhood of $\overline x$. Supposing the contrary, there must exist a $\bar \delta>0$ such that $B(\overline{x},\bar \delta)\subseteq S$. But this is a contradiction to \eqref{s} and, consequently, $S$ fails to be a neighborhood of $\overline{x}$.
\end{example}

\begin{example}\label{extreiman1}
The following examples shows that in the second assertion of Lemma \ref{treiman} one cannot renounce at the hypotheses that $f$ is calm at $\overline{x}$. Indeed, take $S$ a sponge around $\overline{x} \in X$, which is not a neighborhood of $\bar x$ and define $f :X \rightarrow \R$ as being
$$f(x)= \left \{\begin{array}{ll}
                  0, & \ \mbox{if} \ x \in S, \\
                  -1, & \ \mbox{otherwise}.
                  \end{array}\right.$$
Then, for all $\varepsilon \geq 0$, $0$ is an $H_\varepsilon$-subgradient of $f$ at $\bar x$, but $f$ is not calm at $\bar x$ and, consequently,
$0 \notin \partial^-_{\varepsilon} f(\overline{x})$.
\end{example}

\begin{example}\label{extreiman2}
Both assertions of Lemma \ref{treiman} have been given by Treiman in \cite{Treiman} without proof for $f$ a lower semicontinuous function on $X$ and without assuming for (ii) that $f$ is calm at $\bar x$. The following example, which has been kindly provided to us by Jean-Paul Penot, shows that even for lower semicontinuous functions one cannot renounce at the calmness hypotheses in order to get the desired conclusion. Let $X$ be an infinite dimensional Banach space and a sequence of elements $\{e_n\}_{n \geq 1}$ on the unit sphere of $X$ such that $\|e_n - e_m\| > 1/2$ for all $n, m \geq 1$, $n \neq m$. Define $f : X \rightarrow \R$ as being $f(x) = -1/{2^n}$ when $n \geq 1$ is such that $x = 1/{4^n} e_n$ and $f(x) = 0$, otherwise. The function $f$ is lower semicontinuous and it fulfills $f(x) \geq f(0)$ for all $x \in X \setminus \bigcup_{n \geq 1} \left \{1/{4^n} e_n \right \}$. Since $X \setminus \bigcup_{n \geq 1} \left \{1/{4^n} e_n \right \}$ is a sponge around $0$, for all $\varepsilon \geq 0$, $0$ is an $H_\varepsilon$-subgradient of $f$ at $\bar x$. On the other hand,  as $f$ fails to be calm at $0$, $0$ cannot be a Dini-Hadamard $\varepsilon$-subgradient of $f$ at $0$.
\end{example}

Next we provide a variational description of the Dini-Hadamard $\varepsilon$-subdifferential similar to the one that exists for the Fr\' echet $\varepsilon$-subdifferential, but by replacing neighborhoods with sponges.

\begin{theorem}\label{sponge-subdifferential}
Let $f:X \ra \R \cup \{+\infty\}$ be an arbitrary function and $\overline x \in \dom f$. Then for all $\varepsilon \geq 0$ one has
\begin{equation}\label{sponge-dini}
\begin{array}{rl}
x^* \in \partial^-_\varepsilon f(\overline{x}) \Leftrightarrow & f \ \mbox{is calm at} \ \overline{x} \ \mbox{and} \ \forall \alpha>0 \ \mbox{there exists} \ S \ \mbox{a sponge around} \ \overline{x} \ \mbox{such that}\\
&\forall x \in S \ f(x)-f(\overline{x}) \geq \langle x^*, x-\overline{x} \rangle - (\alpha + \varepsilon) \|x-\overline{x}\|.
\end{array}
\end{equation}
\end{theorem}

\begin{proof} Consider an $\varepsilon \geq 0$ fixed.

``$\subseteq$'' Let $x^* \in \partial^-_\varepsilon f(\overline{x})$ and $\alpha>0$. Using Lemma \ref{treiman} and Remark \ref{sponge} one obtains the existence of a sponge $S$ around $\overline{x}$ such that for all $x \in S$
$$f(x)-f(\overline{x}) \geq \langle x^*, x-\overline{x} \rangle - (\alpha + \varepsilon) \|x-\overline{x}\|.$$
Further, from $x^* \in \partial^-_\varepsilon f(\overline{x})$ it follows that $d^-f(\overline{x};0) = 0$ and this provides the desired inclusion.

``$\supseteq$'' For the reverse inclusion assume that $f$ is calm at $\overline{x}$ and consider an arbitrary element $x^*$ fulfilling the property in the right-hand side of \eqref{sponge-dini}. We have to show that
\begin{eqnarray}\label{ineq}
d^-f(\overline{x};h) \geq \langle x^*,h \rangle - \varepsilon \|h\| \ \forall h \in X.
\end{eqnarray}
Let $h \in X \setminus \{0\}$ be fixed. For all $k \in \N$, by taking $\alpha_k:=\frac{1}{k}$, there exists $S_k$ a sponge around $\overline{x}$ such that for all $x \in S_k$
$$f(x)-f(\overline{x}) \geq \langle x^*, x-\overline{x} \rangle - \left(\frac{1}{k} + \varepsilon\right) \|x-\overline{x}\|.$$
Thus, for all $k \in \N$ there exist $\lambda_k>0$ and $\delta_k >0$ such that for all $t \in (0,\lambda_k)$ and all $u \in B(h,\delta_k)$ one has $\overline{x}+tu \in S_k$ and
$$f(\overline{x}+tu)-f(\overline{x}) \geq \langle x^*,tu \rangle - \left(\frac{1}{k}+\varepsilon\right)\|tu\|,$$
which imply in turn that
$$d^-f(\overline{x};h) = \liminf_{\substack{u \ra h \\ t \downarrow 0}} \frac{f(\overline{x} + tu) - f(\overline{x})}{t} \geq$$
$$\liminf_{\substack{u \ra h \\ t \downarrow 0}} \left[ \langle x^*,u \rangle - \left( \frac{1}{k}+\varepsilon \right)\|u\| \right] = \langle x^*,h \rangle - \left( \frac{1}{k}+\varepsilon \right)\|h\|.$$
Passing now to the limit as $k \ra +\infty$, we finally obtain that
$$d^-f(\overline{x};h) \geq \langle x^*,h \rangle - \varepsilon \|h\|.$$
Noticing that, due to the calmness of $f$ at $\overline{x}$, the above inequality is valid also in case $h=0$, the desired conclusion follows.
\end{proof}

\begin{remark}
By making use of Theorem \ref{sponge-subdifferential} one can easily prove that for $\varepsilon \geq 0$ the Dini-Hadamard $\varepsilon$-subdifferential of $f:X \ra \R \cup \{+\infty\}$ at $\overline x \in \dom f$ can be also characterized at follows
\begin{equation*}\label{dini-2}
\begin{array}{rl}
x^* \in \partial^-_\varepsilon f(\overline{x}) \Leftrightarrow& f \ \mbox{is calm at} \ \overline{x} \ \mbox{and} \ \forall \alpha>0 \ \forall u \in S_X \ \exists  \delta >0 \ \mbox{such that}\\
&\forall s \in (0,\delta) \ \forall v \in B(u, \delta) \ \mbox{for} \ x:=\overline{x}+sv \ \mbox{one has}\\
&f(x)-f(\overline{x}) \geq \langle x^*, x-\overline{x} \rangle - (\alpha + \varepsilon) \|x-\overline{x}\|.
\end{array}
\end{equation*}
\end{remark}

Further, let $f, g:X \ra \R \cup \{+\infty\}$ be two arbitrary functions. By using Theorem \ref{sponge-subdifferential} and the fact that the intersection of two sponges around the same point is a sponge around that point, one can prove that for all $\varepsilon, \eta \geq 0$ and all $x \in \dom f \cap \dom g$
\begin{eqnarray}\label{alpha-beta}
\partial^-_\varepsilon f(x)+\partial^-_\eta g(x) \subseteq \partial^-_{\varepsilon+\eta} (f+g)(x).
\end{eqnarray}
From the conventions made for the Dini-Hadamard $\varepsilon$-subdifferential it follows that \eqref{alpha-beta} is in fact true for all $x \in X$.

In what follows we give via \eqref{alpha-beta} a formula for the \textit{difference} of two functions. To this end we need to introduce the notion of \emph{star-difference} of two sets. For $A, B \subseteq X$ the \emph{star-difference} of $A$ and $B$ is defined as
$$A\frac{*}{}B:=\{x \in X: \ x + B \subseteq A \}=\bigcap_{\substack{b \in B}}\{A-b\}.$$
This notion has been introduced by Pontrjagin in \cite{Pontrjagin} in the context of linear differential games and has found resonance in different theoretical and practical investigations in the field of nonsmooth analysis (see, for instance, \cite{Amahroq, Aubin, Caprari, Gautier, Hiriart-Urruty, Martinez-Legaz, MNY, Pshenichnii}).

When dealing with the difference of two functions $g, h: X \ra \R \cup \{+\infty\}$ we assume throughout this paper that $\dom g \subseteq \dom h$. This guarantees that the function $f=g-h : X \rightarrow \R \cup \{+\infty\}$ is well-defined. Moreover, one can easily verify that $g=f+h$ and $\dom f = \dom g$ and, consequently, by making use of \eqref{alpha-beta}, we get the following result.

\begin{proposition}\label{fg-arbitrary}
Let $g, h: X \ra \R \cup \{+\infty\}$  be given functions with $\dom g \subseteq \dom h$ and $f:=g-h$. Then for all $\varepsilon, \eta \geq 0$ and all $x \in X$ one has
\begin{eqnarray}\label{star-difference}
\partial^-_\varepsilon f(x) \subseteq \partial^-_{\varepsilon + \eta} g(x) \frac{*}{} \partial^-_\eta h(x).
\end{eqnarray}
\end{proposition}

\begin{remark}\label{remarkfg} (a) If for $\eta \geq 0$ and $x \in X$ the set $\partial^-_\eta h(x)$ is nonempty, then $\partial^-_{\varepsilon + \eta} g(x) \frac{*}{} \partial^-_\eta h(x)$ $\subseteq \partial^-_{\varepsilon + \eta} g(x) - \partial^-_\eta h(x)$ for all $\varepsilon \geq 0$.

(b) If $\overline{x} \in \dom f$ is a local minimizer of the function $f:=g-h$, then
$$0 \in \partial^-g(\overline{x}) \frac{*}{} \partial^-h(\overline{x})$$
or, equivalently,
$$\partial^-h(\overline{x}) \subseteq \partial^-g(\overline{x}).$$

(c) Similar characterizations for the difference of two functions to the one in Proposition \ref{fg-arbitrary} have been given in \cite{Amahroq} by means of the Fr\'echet subdifferential and in \cite{MNY} by means of the \textit{Mordukhovich (basic/limiting) subdifferential} (see \cite{Morduchovich1, Morduchovich2}).
\end{remark}

\section{The difference of two directionally approximately starshaped functions}

In this section we show first of all that for some particular classes of functions one gets equality in \eqref{star-difference}. After that, we employ these investigations to the formulation of optimality conditions via the Dini-Hadamard subdifferential for cone-constrained optimization problems having the difference of two functions as objective.

We start by presenting some generalized convexity notions for functions.
\begin{definition}\label{ac}
Let $f:X \ra \R \cup \{+\infty\}$ be a given function and $\overline{x} \in \dom f$. The function $f$ is said to be

(i) \emph{approximately convex} at $\overline{x}$, if for any $\varepsilon >0$ there exists $\delta >0$ such that for every $x,y \in B(\overline{x},\delta)$ and every $t \in [0,1]$ one has
\begin{eqnarray}\label{def-ac}
f((1-t)y+tx) \leq (1-t)f(y)+tf(x)+\varepsilon t(1-t) \|x-y\|.
\end{eqnarray}

(ii) \emph{approximately starshaped} at $\overline{x}$, if for any $\varepsilon >0$ there exists $\delta >0$ such that for every $x \in B(\overline{x},\delta)$ and every $t \in [0,1]$ one has
\begin{eqnarray}\label{def-as}
f((1-t)\overline{x}+tx) \leq (1-t)f(\overline{x})+tf(x)+\varepsilon t(1-t) \|x-\overline{x}\|.
\end{eqnarray}

(iii) \emph{directionally approximately starshaped} at $\overline{x}$, if for any $\varepsilon >0$ and any $u \in S_X$ there exists $\delta >0$ such that for every $s \in (0, \delta)$, every $v \in B(u, \delta)$ and every $t \in [0,1]$, when $x:=\overline{x}+sv$, one has
\begin{eqnarray}\label{def-das}
f((1-t)\overline{x}+tx) \leq (1-t)f(\overline{x})+tf(x)+\varepsilon t(1-t) \|x-\overline{x}\|.
\end{eqnarray}
\end{definition}
The approximately convex functions have been introduced in \cite{Ngai-ac} (see also \cite{Aussel, Penot-ac}), while the approximately starshaped  and the directionally approximately starshaped ones have been object of study in \cite{Penot-as-das}.

\begin{remark}\label{remac}
The set of approximately convex functions at a given point $\overline{x} \in X$ is a convex cone containing the functions which are strictly differentiable at $\overline{x}$, being stable under finite suprema (see \cite{Ngai-ac}). An example of an approximately convex function at every $x \in \R$, which is not convex, is
$x \mapsto |x| - x^2$.
\end{remark}

\begin{remark}\label{remacas}
One can easily see that if $f$ is approximately convex at $\overline{x}$, then it is approximately starshaped at $\overline{x}$, too. Nevertheless, the reverse implication does not hold. The following example in this sense has been inspired by \cite[Example 6.10]{Penot-as-das}. We define $f: \R \rightarrow \R$ as follows: $f(0) = 0$, $f(x) = 1/(2n+1) (x-1/(2n)) + 1/(4n^2)$, for $x \in [1/(2n+1), 1/(2n)]$, $n \geq 1$, $f(x) = 1/(2n)x$, for $x \in [1/(2n), 1/(2n-1))$, $n \geq 1$, $f(x) = +\infty$, for $x \geq 1$, while for $x  < 0$ we take $f(x) = f(-x)$.
Then $f$ is approximately starshaped at $0$, but not approximately convex at $0$.
\end{remark}

\begin{remark}\label{remasdas}
By a straightforward calculation one can show that if $f$ is approximately starshaped at $\overline{x}$, then it is directionally approximately starshaped at $\overline{x}$, too. In order to give an example for the failure of the reverse implication we first characterize the class of directionally approximately starshaped functions by means of sponges. A direct consequence of Proposition \ref{das} will be the fact that, in finite dimensional spaces, the two classes of functions coincide (see also \cite{Penot-as-das}).
\end{remark}

\begin{proposition}\label{das}
Let $f:X \ra \R \cup \{+\infty\}$ be a given function and $\overline{x} \in \dom f$. Then $f$ is directionally approximately starshaped at $\overline{x}$ if and only if for any $\varepsilon >0$ there exists a sponge $S$ around $\overline{x}$ such that for every $x \in S$ and every $t \in [0,1]$ one has
\begin{eqnarray}\label{f-das}
f((1-t)\overline{x}+tx) \leq (1-t)f(\overline{x})+tf(x)+\varepsilon t(1-t) \|x-\overline{x}\|.
\end{eqnarray}
\end{proposition}

\begin{proof} As the sufficiency follows directly from the definition of the sponge, we only prove the necessity and assume that $f$ is directionally approximately starshaped at $\overline{x}$. Let be $\varepsilon >0$ fixed. We show that the set
$$S:=\{x \in X: \  f((1-t)\overline{x}+tx) \leq (1-t)f(\overline{x})+tf(x)+\varepsilon t(1-t) \|x-\overline{x}\|, \ \forall t \in [0,1]\}$$
is a sponge around $\overline{x}$. To this end, we take an arbitrary $h \in X \setminus \{0\}$. Since $\frac{1}{\|h\|}h \in S_X$, there exists $\delta >0$ such that for every $s \in (0, \delta)$, every $v \in B\left(\frac{1}{\|h\|}h, \delta\right)$ and every $t \in [0,1]$, when $x:=\overline{x}+sv$, it holds
$$f((1-t)\overline{x}+tx) \leq (1-t)f(\overline{x})+tf(x)+\varepsilon t(1-t) \|x-\overline{x}\|.$$
Let be now $\delta_1:=\frac{\delta}{2\|h\|}>0$ and $\delta_2:=\delta \|h\|>0$. For every $s_1 \in (0, \delta_1]$ and every $v_1 \in B(h, \delta_2)$ one has $s_1 \|h\| \in (0,\delta)$ and $\frac{1} {\|h\|}v_1 \in B\left(\frac{1}{\|h\|}h, \delta\right)$ and, in this way, for every $t \in [0,1]$ and $x:=\overline{x}+ s_1v_1 = \overline{x} + s_1 \|h\|\frac{1}{\|h\|}v_1$ it holds
$$f((1-t)\overline{x}+tx) \leq (1-t)f(\overline{x})+tf(x)+\varepsilon t(1-t) \|x-\overline{x}\|.$$
Consequently, $\overline{x} + [0, \delta_1]B(h,\delta_2) \subseteq S$, which means that $S$ is a sponge around $\overline{x}$ and this leads to the desired conclusion.
\end{proof}

We come now to the announced example of a function which is directionally approximately starshaped at a point, but fails to be approximately starshaped at that point.

\begin{example}\label{das-not-as}
Let $\overline{x} \in X$, $S \subseteq X$ be a sponge around $\overline{x}$, which is not a neighborhood of $\bar x$ (see, for instance, Example \ref{sponge-not-neighborhood}), and the function $f : X \rightarrow \R$,
$$f(x)= \left \{\begin{array}{ll}
                  0, & \ \mbox{if} \ x \in S, \\
                  -\|x-\overline{x}\|, & \ \mbox{otherwise}.
                  \end{array}\right.$$
The function $f$ is directionally approximately starshaped at $\overline{x}$, as for every $\varepsilon > 0$ the inequality in \eqref{f-das} is fulfilled for all $x$ in the given sponge $S$ and all $t \in [0,1]$. We show that $f$ is not approximately starshaped at $\overline{x}$ by assuming the contrary. This means that for $\varepsilon=1/2$ there exists $\delta > 0$ such that for every $x \in B(\overline{x},\delta)$ and every $t \in [0,1]$ one has
\begin{eqnarray}\label{eqexample}
f((1-t)\overline{x}+tx) \leq (1-t)f(\overline{x})+tf(x)+\frac{1}{2} t(1-t) \|x-\overline{x}\|.
\end{eqnarray}
Since $S$ is not a neighborhood of $\bar x$, there exists $x \in B(\overline{x}, \delta) \setminus S$. As $S$ is a sponge around $\overline{x}$ and $x \neq \overline{x}$, there exists $\lambda \in (0,1)$ such that $\overline{x} + [0,\lambda](x-\overline{x}) \in S$. But from \eqref{eqexample} it follows that for all
$t \in (0,\lambda]$
$$t\|x-\overline{x}\| \leq \frac{1}{2}t(1-t)\|x-\overline{x}\|.$$
Dividing now by $t$ and passing to the limit as $t \downarrow 0$, we obtain $\|x-\overline{x}\| \leq 0$ and hence $x=\overline{x}$, which is impossible.
Consequently, $f$ is not approximately starshaped at $\overline{x}$.
\end{example}

\begin{remark}\label{diffFrDi} According to \cite[Lemma 26]{Penot-as-das}, if $f$ is approximately starshaped at $\overline{x} \in \dom f$, then one has
$\partial^F f(\overline{x})=\partial^- f(\overline{x})$. From the previous example it follows that this is no longer true if the function is (only)
directionally approximately starshaped at $\overline{x}$. Indeed, for the function in Example \ref{das-not-as} one gets, since $f$ is calm at $\overline{x}$, via Theorem \ref{sponge-subdifferential}, that
$0 \in \partial^- f(\overline{x})$, while, by employing the characterization \eqref{charFr}, it follows that $0 \notin \partial^F f(\overline{x})$.
\end{remark}

We state the following result on directionally approximately starshaped functions without proof, since it is a direct consequence of \cite[Lemma 27]{Penot-as-das}, whereby the sponge can be constructed like in \eqref{eq7}.

\begin{lemma}\label{das-lemma}
Let the $f:X \ra \R \cup \{+\infty\}$ be directionally approximately starshaped at $\overline{x} \in \dom f$. Then for every $\alpha>0$ and every $\varepsilon \geq 0$ there exists a sponge $S$ around $\overline{x}$ such that for every $x \in S$ one has
\begin{eqnarray}
f(x)-f(\overline{x}) \geq \langle \overline{x}^*,x-\overline{x} \rangle - (\alpha + \varepsilon) \|x-\overline{x}\| \ \forall \overline{x}^* \in \partial^-_\varepsilon f(\overline{x})\\ \label{x-lemma}
f(\overline{x})-f(x) \geq \langle x^*,\overline{x}-x \rangle - (\alpha + \varepsilon) \|x-\overline{x}\| \ \forall x^* \in \partial^-_\varepsilon f(x). \label{xbar-lemma}
\end{eqnarray}
\end{lemma}

The notion we introduce next extends the gap-continuity introduced by Penot in \cite{Penot-gap-continuity}.

\begin{definition}\label{gap-continuity}
\rm A multivalued mapping $F:X \rightrightarrows Y$ between a topological space $X$ and a metric space $Y$ is said to be \emph{spongiously gap-continuous} at $\overline{x} \in X$ if for any $\varepsilon>0$ there exists a sponge $S$ around $\overline{x}$ such that for every $x \in S$
$$\gap(F(\overline{x}), F(x)) < \varepsilon,$$
where for two subsets $A$ and $B$ of $Y$
$$\gap(A,B):=\inf\{d(a,b): a \in A, b \in B\},$$
with the convention that if one of the sets is empty, then $\gap(A,B):=+\infty$.
\end{definition}

When defining a gap-continuous mapping one only has to replace in the above definition the sponge $S$ around $\overline{x}$ with a neighborhood of $\overline{x}$ in $X$. Therefore, every gap-continuous mapping at a point is spongiously gap-continuous at that point. Every multivalued mapping which is either lower semicontinuous or upper semicontinuous at a given point is gap-continuous (see \cite{Penot-gap-continuity}) and, consequently, spongiously gap-continuous at that point. On the other hand, it can be shown that $F:X \rightrightarrows Y$ is spongiously gap-continuous at $\overline{x}$ if and only if for any $\varepsilon>0$ there exists a sponge $S$ around $\overline{x}$ such that for every $x \in S$ one has
\begin{eqnarray}\label{sp-gap-cont}
F(x) \cap (F(\overline{x})+\varepsilon B_Y) \neq \emptyset.
\end{eqnarray}
When $\overline{x} \in X$ and $S \subseteq X$ is a sponge around $\overline{x}$, which is not a neighborhood of $\overline x$, then $F:X \rightrightarrows \R$ defined by $F(x) = \{0\}$ for $x \in S$ and $F(x) = \emptyset$, otherwise, is not gap-continuous, but spongiously gap-continuous at $\overline{x}$.

\begin{proposition}\label{gap-property} Let $F,G:X \rightrightarrows Y$ be two multivalued mappings. If $F$ is spongiously gap-continuous at $\overline{x} \in X$ and there exists a sponge $S$ around $\overline{x}$ such that $F(x) \subseteq G(x)$ for all $x \in S$, then $G$ is spongiously gap-continuous at $\overline{x}$.
\end{proposition}

\begin{proof} Let be $\varepsilon>0$. Using Definition \ref{gap-continuity} we get a sponge $T$ around $\overline{x}$ such that for every $x \in T$
$$\inf \{d(\overline{x}^*,x^*): \ \overline{x}^* \in F(\overline{x}), x^* \in F(x)\} < \varepsilon.$$
The set $S \cap T$ is a sponge around $\overline{x}$, too, and for all $x \in S \cap T$ we have
$$\inf\{d(\overline{x}^*,x^*): \ \overline{x}^* \in G(\overline{x}), x^* \in G(x)\} \leq \inf\{d(\overline{x}^*,x^*): \ \overline{x}^* \in F(\overline{x}), x^* \in F(x)\} < \varepsilon,$$
which concludes the proof.
\end{proof}

\begin{remark}\label{remspgapct}
For $f:X \ra \R \cup \{+\infty\}$ and $\overline{x} \in \dom f$, via the monotonicity property of the Dini-Hadamard $\varepsilon$-subdifferential (cf. Remark \ref{remarca-1}), it follows that $\partial^-_\eta f$ is spongiously gap-continuous at $\overline{x}$ for all $\eta > 0$, whenever $\partial^- f$ is spongiously gap-continuous at $\overline{x}$.
\end{remark}

The following result gives a first refinement of the statement in Proposition \ref{fg-arbitrary} in case $\varepsilon=\eta=0$.

\begin{theorem}\label{T1}
Let $g, h: X \ra \R \cup \{+\infty\}$  be two directionally approximately starshaped functions at $\overline{x} \in \dom g \subseteq \dom h$  such that $\partial^-h$ is spongiously gap-continuous at $\overline{x}$ and $f:=g-h$ is calm at $\overline{x}$. Then it holds
\begin{eqnarray}\label{diff-das}
\partial^- f(\overline{x})=\partial^- g(\overline{x}) \frac{*}{} \partial^- h(\overline{x}).
\end{eqnarray}
\end{theorem}

\begin{proof} In view of Proposition \ref{fg-arbitrary} we only have to prove that $\partial^- g(\overline{x}) \frac{*}{} \partial^- h(\overline{x}) \subseteq \partial^- f(\overline{x})$. To this end take an arbitrary $\overline{x}^* \in \partial^- g(\overline{x}) \frac{*}{} \partial^- h(\overline{x})$ and fix $\alpha >0$.
Using Lemma \ref{das-lemma} and the fact that $\partial^- h(\overline{x})$ is spongiously gap-continuous at $\overline x$, one gets a sponge $S$ around $\overline{x}$ such that for all $x \in S$
\begin{eqnarray}
\exists z^*_x \in \partial^- h(x), \ \exists \overline{z}^*_x \in \partial^- h(\overline{x}) \ \mbox{such that} \ \|z^*_x-\overline{z}^*_x\| \leq \frac{\alpha}{3} \label{1},\\
\langle {u}^*,x-\overline{x} \rangle \leq g(x)-g(\overline{x}) + \frac{\alpha}{3} \|x-\overline{x}\| \ \forall {u}^* \in \partial^- g(\overline{x}) \label{2}
\end{eqnarray}
and
\begin{equation}\label{3}
\langle v^*,\overline{x}-x \rangle \leq h(\overline{x})-h(x) + \frac{\alpha}{3} \|x-\overline{x}\| \ \forall {v}^* \in \partial^- h(x).
\end{equation}
Let be now $x \in S$ fixed. Taking $u^*:=\overline{x}^* + \overline{z}^*_x \in \partial^- g(\overline{x})$ in \eqref{2} and $v^*:=z^*_x \in \partial^- h(x)$ in \eqref{3} and adding the two inequalities, we obtain
$$\langle \overline{x}^*+\overline{z}^*_x-z^*_x, x-\overline{x} \rangle \leq f(x)-f(\overline{x})+\frac{2\alpha}{3} \|x-\overline{x}\|$$
$$\Leftrightarrow \langle \overline{x}^*, x-\overline{x} \rangle \leq f(x)-f(\overline{x}) + \langle z^*_x - \overline{z}^*_x, x-\overline{x} \rangle + \frac{2\alpha}{3} \|x-\overline{x}\|$$
and from here, via \eqref{1},
$$\langle \overline{x}^*, x-\overline{x} \rangle \leq f(x)-f(\overline{x}) + \alpha \|x-\overline{x}\|.$$
By Theorem \ref{sponge-subdifferential} it follows that $\overline{x}^* \in \partial^- f(\overline{x})$, which concludes the proof.
\end{proof}

\begin{remark}\label{remT1}
(a) In the hypotheses of Theorem \ref{T1} one has that $0 \in \partial^- f(\overline{x})$ if and only if $\partial^- h(\overline{x}) \subseteq \partial^- g(\overline{x})$.

(b) If $g, h: X \ra \R \cup \{+\infty\}$  are convex functions with $\dom g \subseteq \dom h$, $\partial^-h$ is spongiously gap-continuous at $\overline{x} \in \dom g$ and $f:=g-h$ is calm at $\overline{x}$, then
$$\partial^- f(\overline{x})=\partial g(\overline{x}) \frac{*}{} \partial h(\overline{x}).$$
\end{remark}

Theorem \ref{T1} is the main ingredient for the proof of the following result.

\begin{theorem}\label{T2}
Let $g, h: X \ra \R \cup \{+\infty\}$  be two directionally approximately starshaped functions at $\overline{x} \in \dom g \subseteq \dom h$ and $f:=g-h$ is calm at $\overline{x}$. If for some $\eta \geq 0$ the multivalued mapping $\partial^-_\eta h$ is spongiously gap-continuous at $\overline{x}$, then for all $\varepsilon \geq 0$ it holds
\begin{eqnarray}\label{epsilon-star-difference}
\partial^-_\varepsilon f(\overline{x})=\partial^-_{\varepsilon+\eta} g(\overline{x}) \frac{*}{} \partial^-_\eta h(\overline{x}).
\end{eqnarray}
\end{theorem}

\begin{proof} Let $\eta \geq 0$ be such that  $\partial^-_\eta h$ is spongiously gap-continuous at $\overline{x}$ and $\varepsilon \geq 0$ be fixed. In view of Proposition \ref{fg-arbitrary} we only have to prove that $\partial^-_{\varepsilon+\eta} g(\overline{x}) \frac{*}{} \partial^-_\eta h(\overline{x}) \subseteq \partial^-_\varepsilon f(\overline{x})$. To this end we consider
an arbitrary $x^*$ in the set on the left-hand side of the inclusion above. Taking $g_{\varepsilon+\eta}, h_\eta : X \rightarrow \R \cup \{+\infty\}$ as being $g_{\varepsilon+\eta}(x) := g(x) + (\varepsilon + \eta) \|x-\overline{x}\|$ and $h_\eta(x) := h(x) + \eta \|x-\overline{x}\|$, respectively, via Lemma \ref{dini fepsilon}, it holds
$x^* \in \partial^- g_{\varepsilon+\eta}(\overline{x}) \frac{*}{} \partial^-h_\eta(\overline{x})$. The functions $g_{\varepsilon+\eta}$ and $h_\eta$ are both directionally approximately starshaped at $\overline{x}$, while $\partial^-h_\eta$ is spongiously gap-continuous at $\overline{x}$. Since $f_\varepsilon(x):= f(x) +\varepsilon \|x-\overline{x}\|=
g_{\varepsilon+\eta}(x) - h_\eta(x)$ for $x \in X$ is calm at $\overline{x}$, by Theorem \ref{T1} and taking again into account Lemma \ref{dini fepsilon}, one obtains
$$\partial^-_{\varepsilon+\eta} g(\overline{x}) \frac{*}{} \partial^-_\eta h(\overline{x}) = \partial^- g_{\varepsilon+\eta}(\overline{x}) \frac{*}{} \partial^- h_\eta(\overline{x})
= \partial^-f_\varepsilon(\overline{x}) = \partial^-_\varepsilon f(\overline{x})$$
and hence the desired conclusion.
\end{proof}

\begin{remark}\label{remT2} (a) One should notice that, in the hypotheses of Theorem \ref{T2}, for all $\varepsilon \geq 0$ it holds
\begin{eqnarray}\label{4}
\partial^-_\varepsilon f(\overline{x})=\bigcap_{\substack{\mu \geq 0}}\left (\partial^-_{\varepsilon+\mu} g(\overline{x}) \frac{*}{} \partial^-_\mu h(\overline{x})\right).
\end{eqnarray}
(b) As pointed out in Remark \ref{remspgapct}, in order to guarantee that $\partial^-_\eta h$ is spongiously gap-continuous at $\overline{x}$ for a given $\eta \geq 0$, it is enough to assume that $\partial^- h$ is spongiously gap-continuous at $\overline{x}$.
\end{remark}

\begin{remark}\label{remT1T2} For similar results to Theorem \ref{T1} and Theorem \ref{T2} expressed by means of the Fr\' echet subdifferential we refer to
\cite[Theorem 1]{Amahroq} and \cite[Theorem 3]{Amahroq}, respectively. There, the functions $g$ and $h$ are assumed to be approximately starshaped at $\overline{x}$ and  the spongious gap-continuity is replaced by the gap-continuity. The fact that there exist directionally approximately starshaped functions which are not approximately starshaped fully motivates the necessity of formulating results like Theorem \ref{T1} and Theorem \ref{T2}. Notice that the assumption $\dom g \subseteq \dom h$ seems to be necessary also in \cite{Amahroq}, in order to guarantee that the difference function takes values in $\R \cup \{+\infty\}$, which is the setting considered in the mentioned article, too.
\end{remark}

The following result is a consequence of Theorem \ref{T2} and \cite[Theorem 3]{Amahroq} (see also Lemma \ref{dini fepsilon}) and describes a situation for which the Fr\' echet and Dini-Hadamard subdifferential of the difference of two approximately starshaped functions at a point coincide (notice the two subdifferentials coincide on approximately starshaped functions, while the difference of two approximately starshaped functions at a point is not necessarily approximately starshaped at that point).

\begin{corollary}\label{corT2T2}
Let $g, h: X \ra \R \cup \{+\infty\}$  be two approximately starshaped functions at $\overline{x} \in \dom g \subseteq \dom h$  with the property that there exists $\eta \geq 0$ such that $\partial_\eta^-h$ is gap-continuous at $\overline{x}$ and $f:=g-h$ is calm at $\overline{x}$. Then for all $\varepsilon \geq 0$ it holds $\partial_\varepsilon^F f(\overline{x}) = \partial_\varepsilon^- f(\overline{x})$.
\end{corollary}

Two further corollaries of Theorem \ref{T2} follow.

\begin{corollary}\label{corT2}
Let $g, h: X \ra \R \cup \{+\infty\}$  be two directionally approximately starshaped functions at $\overline{x} \in \dom g \subseteq \dom h$  such that $\partial^-h$ is spongiously gap-continuous at $\overline{x}$ and $f:=g-h$ is calm at $\overline{x}$. Then the following statements are equivalent:

(i) there exists $\eta \geq 0$ such that $\partial^-_\eta h(\overline{x}) \subseteq \partial^-_\eta g(\overline{x})$;

(ii) $0 \in \partial^- f(\overline{x})$;

(iii) for all $\eta \geq 0$ $\partial^-_\eta h(\overline{x}) \subseteq \partial^-_\eta g(\overline{x})$.
\end{corollary}

\begin{corollary}\label{corT1T2}
Let $g, h: X \ra \R \cup \{+\infty\}$  be two given functions, $\overline{x} \in \dom g \subseteq \dom h$ and $f:=g-h$ be calm at $\overline{x}$. Then the following assertions are true:

(a) If $g$ is convex, $h$ is directionally approximately starshaped at $\overline{x}$ and $\partial^- h$ is spongiously gap-continuous at $\overline{x}$, then for all $\varepsilon \geq 0$ it holds
$$\partial_\varepsilon^- f(\overline{x})=(\partial g(\overline{x})+\varepsilon \overline{B}_{X^*}) \frac{*}{} \partial^- h(\overline{x}).$$

(b) If $g$ is lower semicontinuous, approximately convex at $\overline{x}$, $h$ is directionally approximately starshaped at $\overline{x}$ and $\partial^- h$ is spongiously gap-continuous at $\overline{x}$, then for all $\varepsilon \geq 0$ it holds
$$\partial^-_\varepsilon f(\overline{x})=(\partial^- g(\overline{x})+\varepsilon \overline{B}_{X^*}) \frac{*}{} \partial^- h(\overline{x}).$$
\end{corollary}

\begin{proof} (a) The statement follows via Theorem \ref{T2} and Remark \ref{remarca-1}.

(b) Applying Theorem \ref{T2} for $\eta=0$, we obtain that $\partial_\varepsilon^- f(\overline{x})=\partial_\varepsilon^- g(\overline{x}) \frac{*}{} \partial^- h(\overline{x})$. Since the function $g_\varepsilon : X \ra \R \cup \{+\infty\}$, $g_\varepsilon(x) = g(x) + \varepsilon \|x-\overline{x}\|$ is lower semicontinuous and approximately convex at $\overline{x}$, by \cite[Theorem 3.6]{Ngai-ac}, it follows that (see also Lemma \ref{dini fepsilon} and \eqref{FDG}) $\partial_\varepsilon^- g(\overline{x}) = \partial^- g_\varepsilon(\bar x) = \partial^G g_\varepsilon(\overline{x})$. On the other hand, since $\partial^- g(\overline{x})=\partial^G g(\overline{x})$ (again, by \cite[Theorem 3.6]{Ngai-ac}), from \cite[Corollary 5.6.2]{Ioffe3} one gets $\partial^G g_\varepsilon(\overline{x})=\partial^G g(\overline{x})+\partial^G (\varepsilon \|x-\cdot\|)(\overline{x}) = \partial^- g(\overline{x}) + \varepsilon \overline{B}_{X^*}$ and this leads to the desired conclusion.
\end{proof}

In the final part of the paper we employ the above achievements to the formulation of optimality conditions for a cone-constrained optimization problem having the difference of two function as objective. To this aim we need the following notion.

\begin{definition}\label{def-e-blunt}
Let $C \subseteq X$ be a nonempty set, $f: X \ra \R \cup \{+\infty\}$ be a given function, $\overline{x} \in \dom f \cap C$ and $\varepsilon > 0$. We say that $\overline{x}$ is a \emph{spongiously local $\varepsilon$-blunt minimizer of $f$ on the set $C$} if there exists a sponge $S$ around $\overline{x}$ such that for all $x \in S \cap C$
$$f(x) \geq f(\overline{x})-\varepsilon \|x-\overline{x}\|.$$
In case $C=X$, we simply call $\overline{x}$ a \textit{spongiously local $\varepsilon$-blunt minimizer of $f$}.
\end{definition}

\begin{remark}
It is worth noticing that the above notion generalizes the one of $\varepsilon$-blunt minimizer introduced by Amahroq, Penot and Syam in \cite{Amahroq}. Although in finite dimensional spaces the two notions coincide, this is in general not the case. To see this one only needs to take a look at the Example \ref{das-not-as}. There, $\overline{x}$ is a spongiously local $\varepsilon$-blunt minimizer of $f$ for all $\varepsilon > 0$, but not a local $\varepsilon$-blunt minimizer of $f$ for $\varepsilon \in (0,1)$.
\end{remark}

The following characterization of the Dini-Hadamard subdifferential by means of spongiously local $\varepsilon$-blunt minimizers is a direct consequence of Theorem \ref{sponge-subdifferential}.

\begin{proposition}\label{e-blunt-optimality}
Let $f:X \ra \R \cup \{+\infty\}$ be a given function and $\overline{x} \in \dom f$. Then:
\begin{equation*}
\begin{array}{rl}
0 \in \partial^- f(\overline{x}) \Leftrightarrow & f \ \mbox{is calm at} \ \overline{x} \ \mbox{and} \ \overline{x} \ \mbox{is a spongiously local} \  \varepsilon-\mbox{blunt minimizer of}\\
& f \ \mbox{for all} \ \varepsilon >0.
\end{array}
\end{equation*}
\end{proposition}

Consider now another Banach space $Z$ and $Z^*$ its topological dual space. Let $C \subseteq X$ be a convex and closed set and $K \subseteq Z$ be a nonempty convex and closed \textit{cone} with  $K^*:=\{z^* \in Z^*: \langle z^*,z \rangle \geq 0 \ \mbox{for all} \ z \in K\}$ its \textit{dual cone}. Consider a function $k:X \ra Z$ which is assumed to be \textit{$K$-convex}, meaning that for all $x,y \in X$ and all $t \in [0,1]$, $(1-t)k(x)+tk(y) - k((1-t)x+ty) \in K$, and \textit{$K$-epi closed}, meaning that the \textit{$K$-epigraph} of $k$, $\epi_K k:=\{(x,z) \in X \times Z: z \in k(x)+K\}$, is a closed set. One can notice that when $Z=\R$ and $K=\R_+$ the notion of $K$-epi closedness coincide with the classical lower semicontinuity. For $z^* \in K^*$, by $(z^*k) : X \rightarrow \R$ we denote the function defined by $(z^*k)(x) = \langle z^*, k(x) \rangle$. Further, let $g, h: X \ra \R \cup \{+\infty\}$  be two given functions with $\dom g \subseteq \dom h$ and $f:=g-h$.

The next result provides optimality conditions for the cone-constrained optimization problem
$$\begin{array}{rl}
(\mathcal{P}) & \inf\limits_{x \in \mathcal{A}} f(x).\\
& \mathcal{A}=\{x \in C: k(x) \ \in -K\}
\end{array}
$$

\begin{theorem}\label{optcond}
Let be $\overline{x} \in \inte(\dom g) \cap \mathcal{A}$. Suppose that $g$ is lower semicontinuous and approximately convex at $\overline{x}$, that $f$ is calm at $\overline{x}$ and that $\bigcup_{\substack{\lambda>0}} \lambda (k(C)+K)$ is a closed linear subspace of $Z$. Then the following assertions are true:

(a) If $\overline{x}$ is a spongiously local $\varepsilon$-blunt minimizer of $f$ on $\mathcal{A}$ for all $\varepsilon >0$, then the following relation holds
\begin{eqnarray}\label{condition}
\partial^- h(\overline{x}) \subseteq \partial^- g(\overline{x}) + \bigcup_{\substack{z^* \in K^* \\ (z^*k)(\overline{x}) =0}} \partial((z^*k)+\delta_C)(\overline{x}).
\end{eqnarray}

(b) Viceversa, if $h$ is directionally approximately starshaped at $\overline{x}$, $\partial^-h$ is spongiously gap-continuous at $\overline{x}$ and \eqref{condition} holds, then $\overline{x}$ is a spongiously local $\varepsilon$-blunt minimizer of $f$ on $\mathcal{A}$ for all $\varepsilon >0$.
\end{theorem}

\begin{proof}
We start by noticing that $\mathcal{A}$ is a convex and closed set. Because $\bigcup_{\substack{\lambda>0}} \lambda (k(C)+K)$ is a closed linear subspace of $Z$, via \cite[Theorem 8.6]{Bot-carte}, we have the following representation
\begin{eqnarray*}
N(\mathcal{A}, \overline{x})=\partial(\delta_{\{x \in C: k(x) \in -K\}})(\overline{x})=\bigcup_{\substack{z^* \in K^* \\ (z^*k)(\overline{x})=0}} \partial((z^*k)+\delta_C)(\overline{x})
\end{eqnarray*}
for $N(\mathcal{A}, \overline{x})$, which is the \textit{normal cone} to $\mathcal{A}$ at $\overline{x}$ in the sense of the convex analysis. Therefore, relation \eqref{condition} can be equivalently written as
\begin{eqnarray}\label{condition-1}
\partial^- h(\overline{x}) \subseteq \partial^- g(\overline{x}) + N(\mathcal{A}, \overline{x}).
\end{eqnarray}

On the other hand, taking into consideration that $g$ is lower semicontinuous and approximately convex at $\overline{x} \in \inte(\dom g)$, in view of \cite[Proposition 3.2 and Theorem 3.6]{Ngai-ac}, we have $\partial^- g(\overline{x})=\partial^G g(\overline{x})$, $\partial^-(g+\delta_A)(\overline{x})=\partial^G(g+\delta_A)(\overline{x})$ and $g$ is locally Lipschitz at $\overline{x}$. Further, from \cite[Corollary 5.6.2]{Ioffe3} we get $\partial^G (g+\delta_A)(\overline{x})=\partial^G g(\overline{x})+\partial \delta_\mathcal{A}(\overline{x})=\partial^- g(\overline{x})+N(\mathcal{A}, \overline{x})$, which means that, in the hypotheses of the theorem, relation \eqref{condition} is nothing else than
\begin{eqnarray}\label{g}
\partial^- h(\overline{x}) \subseteq \partial^- (g + \delta_\mathcal{A})(\overline{x}).
\end{eqnarray}

(a) Assuming that $\overline{x} \in \mathcal{A}$ is a spongiously local $\varepsilon$-blunt minimizer of $f$ on $\mathcal{A}$ for all $\varepsilon>0$ means nothing else than that $\overline{x}$ is a spongiously local $\varepsilon$-blunt minimizer of $f+\delta_\mathcal{A}$ for all $\varepsilon>0$ or, equivalently, via Proposition \ref{e-blunt-optimality}, that $0 \in \partial^-(f+\delta_\mathcal{A})(\overline{x}) = \partial^-((g+\delta_\mathcal{A})-h)(\overline{x})$.  From here, by using Proposition \ref{fg-arbitrary}, we obtain relation \eqref{g} and this closes the proof of the first statement.

(b) For proving the second statement we have to notice first that $g+\delta_\mathcal{A}$ is approximately convex at $\overline{x}$. Thus, via Theorem \ref{T1}, one has that $\partial^-((g+\delta_\mathcal{A})-h)(\overline{x}) = \partial^- (g+\delta_\mathcal{A})(\overline{x}) \frac{*}{} \partial^- h(\overline{x})$ and, since \eqref{g} is fulfilled, we have $0 \in \partial^-((g+\delta_\mathcal{A})-h)(\overline{x}) = \partial^-(f+\delta_\mathcal{A})(\overline{x})$. Taking again into account Proposition \ref{e-blunt-optimality} we obtain that $\overline{x}$ is a spongiously local $\varepsilon$-blunt minimizer of $f$ on $\mathcal{A}$ for all $\varepsilon>0$.
\end{proof}

\begin{remark}\label{lastrm}
For a similar result to Theorem \ref{optcond}, given in the particular instance when $K=\{0\}$ and $k(x) = 0$ for all $x \in X$ and by means of the Fr\' echet subdifferential, we refer to \cite[Proposition 6]{Amahroq}. In the second statement of that result the authors ask for $h$ to be approximately starshaped at $\overline{x}$ with $\partial^F h$ gap-continuous at $\overline{x}$ and characterize the local $\varepsilon$-blunt minimizers of $f$ for all $\varepsilon >0$. To this end they make use of some exact subdifferential formulae for the limiting subdifferential, but by providing an incorrect argumentation, since these are valid in Asplund spaces. Nevertheless, the statement in \cite[Proposition 6]{Amahroq} is true in Banach spaces, too, and it can be proven in the lines of the proof of Theorem \ref{optcond}.
\end{remark}

\noindent{\bf Acknowledgements.} The authors are grateful to  E. R. Csetnek for pertinent comments on an earlier draft of the paper.

\end{document}